\documentstyle[12pt]{article} 
\begin{document} 
\newtheorem{Def}{Definition} 
\newtheorem{thm}{Theorem} 
\newtheorem{lem}{Lemma} 
\newtheorem{rem}{Remark} 
\newtheorem{prop}{Proposition} 
\newtheorem{cor}{Corollary} 
\newtheorem{clm}{Claim} 
\newtheorem{step}{Step} 
\newtheorem{sbsn}{Subsection} 
\newtheorem{conj}{Conjection} 
\newtheorem{eg}{Example} 
\title 
{Further results on Liouville type theorems for some conformally 
invariant fully nonlinear equations}  
\author{\ Aobing Li 
\ \ \ \ \& \ \ \ YanYan Li\thanks{Partially supported by 
NSF Grant DMS-0100819.} 
\\ Department of Mathematics\\ Rutgers University\\ 
110 Frelinghuysen Rd.\\ 
Piscataway, NJ 08854 
} 
\date{} 
\maketitle 
\input { amssym.def} 
 
For $n\ge 3$, let ${\cal S}^{n\times n}$ be the set of $n\times n$ 
real symmetric matrices, ${\cal S}_{+}^{n\times n}\subset {\cal 
S}^{n\times n}$ be the set of positive definite matrices, $O(n)$ be 
the set of $n\times n$ real orthogonal matrices. \newline

For a positive $C^2$ function $u$, let 
\[ 
A^u:= -\frac{2}{n-2}u^{  -\frac {n+2}{n-2} } 
\nabla^2u+ \frac{2n}{(n-2)^2}u^ { -\frac {2n}{n-2} } 
\nabla u\otimes\nabla u-\frac{2}{(n-2)^2} u^ { -\frac {2n}{n-2} } 
|\nabla u|^2I, 
\] 
where $I$ is the $n\times n$ identity matrix.

Let $U\subset {\cal S}^{n\times n}$  be an open set satisfying 
\begin{equation} 
O^{-1}UO=U,\qquad\forall\ O\in O(n), 
\label{1} 
\end{equation} 
and 
\begin{equation} 
U\cap \{M+tN\ |\ 0<t<\infty\}\ \mbox{is convex}\qquad \forall\ 
M\in {\cal S}^{n\times n}, N\in {\cal S}^{n\times n}_+. 
\label{2}   
\end{equation}

Let $F\in C^{\infty}(U)$ satisfy 
\begin{equation} 
F(O^{-1}MO)=F(M),\qquad \forall\ 
M\in U, 
\label{3} 
\end{equation} 
\begin{equation} 
\left(F_{ij}(M)\right)>0,\qquad \forall\ M\in U, 
\label{4} 
\end{equation} 
where $F_{ij}(M):=\frac{\partial F}{ \partial M_{ij} }(M)$.\newline

For $n\ge 3$, $-\infty<p\le\frac{n+2}{n-2}$, we consider 
\begin{equation} 
F(A^u)=u^{p-\frac{n+2}{n-2}},\quad A^u\in U,\quad u>0\quad\mbox{on}~\Bbb{R}^n. 
\label{5} 
\end{equation}

Our main theorem is  
\begin{thm} 
For $n\ge 3$, let $U\subset{\cal S}^{n\times n}$ satisfy (\ref{1}), 
(\ref{2}), and let $F\in C^2(U)$ satisfy (\ref{3}), (\ref{4}). Assume 
that $u\in C^2(\Bbb{R}^n)$ is a superharmonic solution of (\ref{5}) 
for some $-\infty<p\le\frac{n+2}{n-2}$. Then either $u\equiv constant$ 
or $p=\frac{n+2}{n-2}$ and for some $\bar x\in\Bbb{R}^n$ and some 
positive constants $a$ and $b$ satisfying $2b^2a^{-2}I\in U$ and 
$F(2b^2a^{-2}I)=1$, 
\begin{equation} 
u(x)\equiv (\frac{a}{1+b^2|x-\bar 
x|^2})^{\frac{n-2}{2}},\quad\forall~x\in\Bbb{R}^n.  
\label{6} 
\end{equation} 
\label{theorem1} 
\end{thm} 
 
\begin{rem} About half a year ago, 
we established a slightly weaker version of  
Theorem \ref{theorem1} for $p<\frac {n+2}{n-2}$, and the proof was different than the one in 
the present paper.  The weaker result requires  
some additional (though minor, e.g.,  $F$ being homogeneous of degree 
$1$ would be enough) assumptions on $(F,U)$. 
\end{rem} 
Theorem~\ref{theorem1} for $p=\frac{n+2}{n-2}$ was established in 
\cite{LL3}, which extends earlier Liouville type theorems for conformally 
invariant equations by Obata (\cite{O}), Gidas, Ni and Nirenberg (\cite{GNN}), 
Caffarelli, Gadis and Spruck (\cite{CGS}), Viaclovsky (\cite{V3} and \cite{V2}), Chang, 
Gursky and Yang (\cite{CGY2} and \cite{CGY3}), and Li and Li (\cite{LL0}, 
\cite{LL}, \cite{LL2} and \cite{LL3}).

The proof of Theorem~\ref{theorem1} for $p=\frac{n+2}{n-2}$ in the 
present paper is simplier than that in our earlier paper \cite{LL3}, though 
the most crucial ideas are the same. Theorem~\ref{theorem1} for 
$-\infty<p<\frac{n+2}{n-2}$ extends the corresponding result of Gidas 
and Spruck in \cite{GS}. The proof of Theorem~\ref{theorem1} for 
$-\infty<p<\frac{n+2}{n-2}$ is essentially the same as our simplified 
proof of Theorem~\ref{theorem1} for $p=\frac{n+2}{n-2}$ in this 
paper. Our proof of Theorem~\ref{theorem1} makes use of the following 
 lemma used in our first proof of Theorem~\ref{theorem1} for 
$p=\frac{n+2}{n-2}$ (see theorem~1 in \cite{LL3}).\newline 
\begin{lem}(\cite{LL3})\ For $n\ge 1$, $R>0$, let $u\in C^2(B_R\setminus\{0\})$ 
satisfying $\Delta u\le 0$ in $B_R\setminus\{0\}$. Assume that there 
exist $w,~v\in C^1(B_R)$ satisfying  
\[ 
w(0)=v(0),\quad\nabla w(0)\neq\nabla v(0), 
\] 
and 
\[ 
u\ge w,\quad u\ge v,\quad\mbox{in}~B_R\setminus\{0\}. 
\] 
Then  
\[ 
\liminf\limits_{x\to 0}u(x)>w(0). 
\] 
\label{lemma0} 
\end{lem} 
 
In fact, the above lemma was stated as lemma 2 in \cite{LL3}  
under additional hypotheses ($w, v\in C^2(B_R)$  
and $\Delta w\le 0$, $\Delta v\le 0$ 
in $B_R$).  However the proof of lemma 2 in \cite{LL3} did not 
use these extra hypotheses.  Indeed, lemma 1 in \cite{LL3} 
was first established and hypothesis (11) there was not used in 
the proof.  So the proof of lemma 2 in \cite{LL3} actually establishes 
Lemma \ref{lemma0} above. 
 
\bigskip

\noindent 
{\bf Proof of Theorem~\ref{theorem1} for $p=\frac{n+2}{n-2}$.}\ Since 
$u$ is a positive superharmonic function, we have, by the maximum 
principle, that 
\[ 
u(x)\ge\frac{\min\limits_{\partial B_1}u}{|x|^{n-2}},\quad\forall~|x|\ge 1. 
\] 
In particular  
\begin{equation} 
\liminf\limits_{|x|\to\infty}(|x|^{n-2}u(x))>0. 
\label{7} 
\end{equation} 
\begin{lem} 
For any $x\in\Bbb{R}^n$, there exists $\lambda_0(x)>0$ such that  
\[ 
u_{x,\lambda}(y):=(\frac{\lambda}{|y-x|})^{n-2} 
u(x+\frac{\lambda^2(y-x)}{|y-x|^2})\le 
u(y),\quad\forall~|y-x|\ge\lambda,~0<\lambda<\lambda_0(x).  
\] 
\label{lemma1} 
\end{lem} 
{\bf Proof of Lemma~\ref{lemma1}.}\ This follows from the proof of 
lemma~2.1 in \cite{LZ}.

\vskip 5pt 
\hfill $\Box$ 
\vskip 5pt 
 
For any $x\in\Bbb{R}^n$, set  
\[ 
\bar\lambda(x):=\sup\{\mu~|~u_{x,\lambda}(y)\le 
u(y),~\forall~|y-x|\ge\lambda,~0<\lambda<\mu\}.  
\] 
Let  
\begin{equation} 
\alpha:=\liminf\limits_{|x|\to\infty}(|x|^{n-2}u(x)). 
\label{8} 
\end{equation} 
Because of (\ref{7}),  
\begin{equation} 
0<\alpha\le\infty. 
\label{8a} 
\end{equation} 
If $\alpha=\infty$, then the moving sphere procedure can never stop 
and therefore $\bar\lambda(x)=\infty$ for any $x\in\Bbb{R}^n$. This 
follows from arguments in \cite{LZ}, \cite{LL} and \cite{LL2}. By the definition of 
$\bar\lambda(x)$ and the fact $\bar\lambda(x)=\infty$, we have,  
\[ 
u_{x,\lambda}(y)\le u(y),\quad\forall~|y-x|\ge\lambda>0. 
\] 
By a calculus lemma (see e.g., lemma~11.2 in \cite{LZ}), $u\equiv 
constant$, and Theorem~\ref{theorem1} for $p=\frac{n+2}{n-2}$ is 
proved in this case (i.e. $\alpha=\infty$). So, from now on, we assume 
 
\begin{equation} 
0<\alpha<\infty. 
\label{9} 
\end{equation} 
By the definition of $\bar \lambda(x)$, 
\[ 
u_{x,\lambda}(y)\le 
u(y),\quad\forall~|y-x|\ge\lambda,~0<\lambda<\bar\lambda (x). 
\] 
Multiplying the above by $|y|^{n-2}$ and sending $|y|\to\infty$, we 
have, 
\[ 
\alpha\ge\lambda^{n-2}u(x),\qquad\forall\ 0<\lambda<\bar\lambda(x). 
\] 
Sending $\lambda\to\bar\lambda(x)$, we have (using (\ref{9})), 
\begin{equation} 
\infty>\alpha\ge\bar\lambda(x)^{n-2}u(x),\quad\forall~x\in\Bbb{R}^n. 
\label{10} 
\end{equation} 
Since the moving sphere procedure stops at $\bar\lambda(x)$, we must 
have, by using the arguments in \cite{LZ}, \cite{LL} and \cite{LL2},  
\begin{equation} 
\liminf\limits_{|y|\to\infty}(u(y)-u_{x,\bar\lambda(x)}(y))|y|^{n-2}=0, 
\label{11} 
\end{equation} 
i.e., 
\begin{equation} 
\alpha=\bar\lambda(x)^{n-2}u(x),\quad\forall~x\in\Bbb{R}^n. 
\label{12} 
\end{equation} 
Let us switch to some more convenient notations. For a Mobius 
transformation $\phi$, we use notation  
\[ 
u_{\phi}:=|J_{\phi}|^{\frac{n-2}{2n}}(u\circ\phi), 
\] 
where $J_{\phi}$ denotes the Jacobian of $\phi$.

For $x\in\Bbb{R}^n$, let 
\[ 
\phi^{(x)}(y):=x+\frac{\bar\lambda(x)^2(y-x)}{|y-x|^2}, 
\]  
we know that 
$u_{\phi(x)}=u_{x,\bar\lambda(x)}$.

Let $\psi(y):=\frac{y}{|y|^2}$, and let 
\[ 
w^{(x)}:=(u_{\phi^{(x)}})_{\psi}=u_{\phi^{(x)}\circ\psi}. 
\]  
For $x\in\Bbb{R}^n$, the only possible singularity for $w^{(x)}$ (on 
$\Bbb{R}^n\cup\{\infty\}$) is $\frac{x}{|x|^2}$. In particular, $y=0$ 
is a regular point of $w^{(x)}$. A direct calculation yields 
\[ 
w^{(x)}(0)=\bar\lambda(x)^{n-2}u(x), 
\] 
and therefore, by (\ref{12}), 
\begin{equation} 
w^{(x)}(0)=\alpha,\quad\forall~x\in\Bbb{R}^n. 
\label{13} 
\end{equation} 
Clearly, $u_{\psi}\in C^2(\Bbb{R}^n\setminus\{0\})$, $\Delta 
u_{\psi}\le 0$ in $\Bbb{R}^n\setminus\{0\}$. We also know that 
\[ 
w^{(x)}(0)=\alpha\quad\forall~x\in\Bbb{R}^n,\qquad\liminf\limits_{y\to 
0}u_{\psi}(y)=\alpha,  
\] 
and, for some $\delta(x)>0$,  
\[ 
w^{(x)}\in C^2(B_{\delta (x)}),\quad\forall~x\in\Bbb{R}^n, 
\] 
\[ 
u_{\psi}\ge 
w^{(x)}\quad\mbox{in}~B_{\delta(x)}\setminus\{0\},\quad\forall~x\in\Bbb{R}^n,  
\] 
\[ 
\Delta w^{(x)}\le 0\quad\mbox{in}~B_{\delta(x)},\quad\forall~x\in\Bbb{R}^n. 
\] 
\begin{lem} 
$\nabla w^{(x)}(0)=\nabla w^{(0)}(0)$, i.e., $\nabla w^{(x)}(0)$ is 
independent of $x\in\Bbb{R}^n$. 
\label{lemma2} 
\end{lem} 
{\bf Proof of Lemma~\ref{lemma2}.}\ This follows from Lemma 
\ref{lemma0}. Indeed, 
for any $x,~\tilde x\in\Bbb{R}^n$, let 
\[ 
v:=w^{(x)},\quad w:=w^{(\tilde x)},\quad u:=u_{\psi}. 
\] 
We know that $w(0)=v(0)$, $u_{\psi}\ge w$ and $u_{\psi}\ge v$ 
near the origin, and we also know that 
 $\liminf\limits_{y\to 0}u_{\psi}(y)=w(0)$, so, 
by Lemma \ref{lemma0},  we must have $\nabla 
v(0)=\nabla w(0)$, i.e., $\nabla w^{(x)}(0)=\nabla w^{(\tilde 
x)}(0)$. Lemma~\ref{lemma2} is established. 
 
\vskip 5pt 
\hfill $\Box$ 
\vskip 5pt

For $x\in\Bbb{R}^n$, 
\begin{eqnarray*} 
w^{(x)}(y)&=&\frac{1}{|y|^{n-2}} 
\Big\{(\frac{\bar\lambda(x)}{|\frac{y}{|y|^2}-x|})^{n-2}u(x+ 
\frac{\bar\lambda(x)^2(\frac{y}{|y|^2}-x)}{|\frac{y}{|y|^2}-x|^2})\Big\}\\ 
&=&(\frac{\bar\lambda(x)}{|\frac{y}{|y|}-|y|x|})^{n-2}u(x+ 
\frac{\bar\lambda(x)^2(y-|y|^2 x)}{|\frac{y}{|y|}-|y|x|^2})\\ 
&=&(\frac{\bar\lambda(x)^2}{1-2x\cdot y+|y|^2 x})^{\frac{n-2}{2}}u(x+ 
\frac{\bar\lambda(x)^2(y-|y|^2x)}{1-2x\cdot y+|y|^2|x|^2}). 
\end{eqnarray*} 
So, for $|y|$ small, 
\[ 
w^{(x)}(y)=\bar\lambda(x)^{n-2}(1+(n-2)x\cdot y)u(x+\bar\lambda(x)^2 y)+O(|y|^2), 
\] 
and, using (\ref{12}), 
\[ 
\nabla 
w^{(x)}(0)=(n-2)\bar\lambda(x)^{n-2}u(x)x+\bar\lambda(x)^n\nabla 
u(x)=(n-2)\alpha x+\alpha^{\frac{n}{n-2}}u(x)^{\frac{n}{2-n}}\nabla u(x). 
\] 
By Lemma~\ref{lemma2}, $\vec V:=\nabla w^{(x)}(0)$ is a constant 
vector in $\Bbb{R}^n$, so we have, 
\[ 
\nabla_x(\frac{n-2}{2}\alpha^{\frac{n}{n-2}}u(x)^{-\frac{2}{n-2}} 
-\frac{(n-2)\alpha}{2}|x|^2+\vec V\cdot x)\equiv 0. 
\] 
Consequently, for some $\bar x\in\Bbb{R}^n$ and $d\in R$,  
\[ 
u(x)^{-\frac{2}{n-2}}\equiv\alpha^{-\frac{2}{n-2}}|x-\bar 
x|^2+d\alpha^{-\frac{2}{n-2}}.  
\] 
Since $u>0$, we must have $d>0$. Thus 
\[ 
u(x)\equiv (\frac{\alpha^{\frac{2}{n-2}}}{d+|x-\bar x|^2})^{\frac{n-2}{2}}. 
\] 
Let $a=\alpha^{\frac{2}{n-2}}d^{-1}$ and $b=d^{-\frac 12}$. Then $u$ 
is of the form (\ref{6}). Clearly $A^u(0)=2b^2a^{-2}I$, so 
$2b^2a^{-2}I\in U$ and $F(2b^2a^{-2}I)=1$. Theorem~\ref{theorem1} in 
the case $p=\frac{n+2}{n-2}$ is established.

\vskip 5pt 
\hfill $\Box$ 
\vskip 5pt 
 
{\bf Proof of Theorem~\ref{theorem1} for 
$-\infty<p<\frac{n+2}{n-2}$.}\   In this case, the equation satisfied by 
$u$ is no longer conformally invariant, but it transforms to our 
advantage when making reflections with respect to spheres, i.e., the 
inequalities have the right direction so that the strong maximum 
principle and the Hopf lemma can still be applied.\newline 
First, we still have (\ref{7}) since this only requires the 
superharmonicity and the positivity of $u$. Lemma~\ref{lemma1} still 
holds since it only uses (\ref{7}) and the $C^1$ regularity of $u$ in 
$\Bbb{R}^n$. For $x\in\Bbb{R}^n$, we still define $\bar\lambda(x)$ in 
the same way. We also define $\alpha$ as in (\ref{8}) and we still 
have (\ref{8a}).\newline 
For $x\in\Bbb{R}^n$, $\lambda>0$, the equation of $u_{x,\lambda}$ now 
takes the form 
\begin{equation} 
F(A^{u_{x,\lambda}}(y))=(\frac{\lambda}{|y-x|})^{(n-2)(\frac{n+2}{n-2}-p)} 
u_{x,\lambda}(y)^{p-\frac{n+2}{n-2}},\quad A^{u_{x,\lambda}}(y)\in 
U,\quad\forall~y\neq x. 
\label{14} 
\end{equation} 
\begin{lem} 
If $\alpha=\infty$, then $\bar\lambda(x)=\infty$ for any $x\in\Bbb{R}^n$. 
\label{lemma3} 
\end{lem} 
{\bf Proof of Lemma~\ref{lemma3}.}\ Suppose the contrary, 
$\bar\lambda(\bar x)<\infty$ for some $\bar x\in\Bbb{R}^n$. Without 
loss of generality, we may assume $\bar x=0$, and we use notations  
\[ 
\bar\lambda:=\bar\lambda(0),\quad u_{\lambda}:=u_{0,\lambda},\quad 
B_{\lambda}:= B_{\lambda}(0). 
\]  
By the definition of $\bar\lambda$,  
\begin{equation} 
u_{\bar\lambda}\le u\quad\mbox{on}~~\Bbb{R}^n\setminus B_{\bar\lambda}. 
\label{15} 
\end{equation} 
By (\ref{14}), 
\begin{equation} 
F(A^{u_{\bar\lambda}})\le u_{\bar\lambda}^{p-\frac{n+2}{n-2}},\quad 
A^{u_{\bar\lambda}}\in U,\quad\mbox{on}~\Bbb{R}^n\setminus B_{\bar\lambda}. 
\label{16} 
\end{equation} 
Recall that $u$ satisfies  
\begin{equation} 
F(A^u)=u^{p-\frac{n+2}{n-2}},\quad A^u\in 
U,\quad\mbox{on}~~\Bbb{R}^n\setminus B_{\bar\lambda}. 
\label{17} 
\end{equation} 
By (\ref{16}) and (\ref{17}), 
\begin{equation} 
F(A^{u_{\bar\lambda}})-F(A^u) 
-(u_{\bar\lambda}^{p-\frac{n+2}{n-2}}-u^{p-\frac{n+2}{n-2}})\le 
0,\quad A^{u_{\bar\lambda}}\in U,~A^u\in 
U,\quad\mbox{on}~~\Bbb{R}^n\setminus B_{\bar\lambda}. 
\label{18} 
\end{equation} 
Since $\alpha=\infty$, we have  
\begin{equation} 
\liminf\limits_{|y|\to\infty}|y|^{n-2}(u-u_{\bar\lambda})(y)>0. 
\label{19} 
\end{equation} 
The inequality in (\ref{18}) goes the right direction. Thus, with 
(\ref{19}), the arguments for $p=\frac{n+2}{n-2}$ work essentially in 
the same way here and we obtain a contradiction by continuing the 
moving sphere procedure a little bit further. This deserves some 
explanations. Because of (\ref{19}), and using arguments in \cite{LL} and 
\cite{LL2}, we only need to show that 
\begin{equation} 
u_{\bar\lambda}(y)< u(y),\quad\forall~|y|>\bar\lambda, 
\label{19a} 
\end{equation} 
and 
\begin{equation} 
\frac{d}{dr}(u-u_{\bar\lambda})|_{\partial B_{\bar\lambda}}>0, 
\label{19b} 
\end{equation} 
where $\frac{d}{dr}$ denotes the differentiation in the outer normal 
direction with repect to $\partial B_{\bar\lambda}$.\newline 
If $u_{\bar\lambda}(\bar y)= u(\bar y)$ for some $|\bar 
y|>\bar\lambda$, then, using (\ref{18}) as in the proof of lemma~2.1 
in \cite{LL}, we know that $u_{\bar\lambda}-u$ satisfies that 
\[ 
L(u_{\bar\lambda}-u)\le 0, 
\] 
where $L=-a_{ij}(x)\partial_{ij}+b_i(x)\partial_i+c(x)$ with 
$(a_{ij})>0$ continuous and $b_i$, $c$ continuous.\newline 
Since $u_{\bar\lambda}-u\le 0$ near $\bar y$, we have, by the strong 
maximum principle, $u_{\bar\lambda}\equiv u$ near $\bar y$. For the 
same reason, $u_{\bar\lambda}(y)\equiv u(y)$ for any 
$|y|\ge\bar\lambda$, violating (\ref{19}). (\ref{19a}) has been 
checked. Estimate (\ref{19b}) can be established in a similar way by 
using the Hopf lemma (see the proof of lemma~2.1 in \cite{LL}). Thus 
Lemma~\ref{lemma3} is established. 
 
\vskip 5pt 
\hfill $\Box$ 
\vskip 5pt

By Lemma~\ref{lemma3} and the usual arguments, we know that if 
$\alpha=\infty$, $u$ must be a constant, and Theorem~\ref{theorem1} 
for $-\infty<p<\frac{n+2}{n-2}$ is also proved in this case.\newline 
>From now on, we always assume (\ref{9}). As before, we obtain 
(\ref{10}). Since the inequality in (\ref{16}) goes the right 
direction, the arguments for $p=\frac{n+2}{n-2}$ (see also the 
arguments in the proof of Lemma~\ref{lemma3}) essentially apply and we 
still have (\ref{11}) and (\ref{12}). The rest of the arguments for 
$p=\frac{n+2}{n-2}$ apply and we have that $u$ is of the form 
(\ref{6}) with some positive constants $a$ and $b$. However, we know 
that, for $u$ of the form (\ref{6}), $A^u\equiv 2b^2a^{-2}I$ and 
$F(A^u)\equiv constant$. This violates (\ref{5}) since 
$u^{p-\frac{n+2}{n-2}}$ is not a constant (recall that 
$p<\frac{n+2}{n-2}$). Theorem~\ref{theorem1} for 
$-\infty<p<\frac{n+2}{n-2}$ is established.  
 
\vskip 5pt 
\hfill $\Box$

\end{document}